\documentclass[runningheads,onecollarge,leqno]{svjour2hack}

\smartqed

\usepackage[USenglish]{babel}
\usepackage{amsmath, amssymb}
\usepackage[algo2e,boxed]{algorithm2e}
\usepackage{algorithmic}

\newcommand{\opt}{{\rm opt}}
\newcommand{\calF} {{{\mathcal F}}}
\newcommand{\calS} {{{\mathcal S}}}
\newcommand{\calP} {{{\mathcal P}}}
\newcommand{\calQ} {{{\mathcal Q}}}
\newcommand{\R}{\mathbb{R}}
\newcommand{\Z}{\mathbb{Z}}

\newcommand*{\transpose}[1]{\ensuremath%
{{#1}^\top}}
\newcommand{\vecone}{\mathbf 1}%
\newcommand{\veczero}{\mathbf 0}%
\newcommand{\conv}{\mbox{conv}}

\newcommand{\argmax}{\mbox{argmax}}
\renewcommand{\vert} {{\rm vert}}
\newcommand{\ignore}[1]{}

\begin{document}

\title{Nonlinear optimization for matroid intersection\\ and extensions} 

\dedication{For Tom Liebling, on the occasion of his retirement}

\author{
Y. Berstein
\and
J. Lee
\and
S. Onn
\and
R. Weismantel
}

\authorrunning{Berstein, Lee, Onn \& Weismantel} % if too long for running head

 \institute{
            Yael Berstein \at
            Technion - Israel Institute of Technology, 32000 Haifa, Israel\\
            \email{yaelber@tx.technion.ac.il}
\and
           Jon Lee \at
           IBM T.J. Watson Research Center, Yorktown Heights, NY 10598, USA\\
           \email{jonlee@us.ibm.com}
\and
           Shmuel Onn \at
           Technion - Israel Institute of Technology, 32000 Haifa, Israel\\
           \email{onn@ie.technion.ac.il}
\and
           Robert Weismantel \at
           Otto-von-Guericke Universit\"at Magdeburg, D-39106 Magdeburg, Germany\\
           \email{weismantel@imo.math.uni-magdeburg.de}
}

\date{}

\maketitle

\begin{abstract}
We address optimization of nonlinear functions of the form $f(Wx)$~,
where $f:\R^d\rightarrow \R$ is a nonlinear function, $W$ is a $d\times n$ matrix,
and feasible $x$ are in some large finite set $\calF$ of integer points
in $\R^n$~. Generally, such problems are intractable, so we obtain positive
algorithmic results by looking at broad natural classes of $f$~, $W$ and $\calF$~.

One of our main motivations is multi-objective discrete
optimization, where $f$ trades off the linear functions given by the
rows of $W$~. Another motivation is that we want to extend as much as possible
the known results about polynomial-time linear optimization over trees,
assignments, matroids, polymatroids, etc. to nonlinear optimization over such structures.

We assume
that the convex hull of $\calF$ is well-described by linear inequalities
(i.e., we have an efficient separation oracle). For example, the set of
characteristic vectors of common bases of a pair of matroids on a common ground set satisfies
this property for $\calF$~.
In this setting, the problem is already known to be
intractable (even for a single matroid), for general $f$ (given by a comparison oracle),
for (i) $d=1$ and binary-encoded $W$~,
and for (ii) $d=n$ and $W=I$~.

Our main results (a few technicalities suppressed):

1- When $\calF$ is well described, $f$ is convex (or even quasiconvex), and $W$ has a fixed number of rows and is unary encoded or
with entries in a fixed set, we give an efficient deterministic algorithm for maximization.

2- When $\calF$ is well described, $f$ is a norm, and binary-encoded $W$ is nonnegative, we give an efficient
deterministic constant-approximation algorithm for maximization.

3- When $\calF$ is well described, $f$ is ``ray concave'' and non-decreasing, and $W$ has a fixed number of rows
and is unary encoded or with entries in a fixed set, we give an efficient deterministic
constant-approximation algorithm for minimization.

4- When $\calF$ is the set of characteristic vectors of common bases of a pair of vectorial matroids on a
common ground set, $f$ is arbitrary, and $W$ has a fixed number of rows and is unary encoded,
we give an efficient randomized algorithm for optimization.
\end{abstract}

\newpage

\section*{Introduction}
Generally, we are considering nonlinear discrete optimization problems of the form
\[
\tag*{$P(\calF,\opt,f,W)$~:}
 \opt~ \left\{f(Wx)~ : ~x\in  \calF\right\}~,
\]
where $\opt\in\{\min,\max\}$~,
$W\in \Z^{d\times n}$ is is a matrix of integers,
function $f:\R^d\rightarrow\R$ is given by a
comparison oracle,
and $\calF$ is a large finite subset of $\Z^n_+$~.

 One motivation for this setting is multiobjective optimization,
 where we  seek an optimal ``balancing'', via
 a nonlinear function $f$~, of $d$ competing linear
 objectives $W_{i\cdot}x$~. Another motivation is the desire to extend
known results about polynomial-time linear optimization
over combinatorial structures to nonlinear optimization
over such structures.

We assume that $\calF$ is ``well described'' --- that is,
 $\calF\subset\calS(n,\beta):=\{ x\in \R_+^n\ :\ \transpose{\vecone}x \le \beta\}~,$
 for some unary encoded $\beta \in \Z_+$~, and that
 we have a separation oracle for $\calP:=\conv(\calF)$~.
 For example, the set of
characteristic vectors of common bases of a pair of matroids on the common ground set
$\{1,2,\ldots,n\}$ satisfies
this property for $\calF$ (with $\beta=n$)~.

 In this setting, the problem is already known to be
intractable (even for a single matroid), for general $f$ (given by a comparison oracle),
for (i) $d=1$ and binary-encoded $W$~,
and for (ii) $d=n$ and $W=I$~.  Therefore, in much of our work, we assume
that the objective vectors $W_{i\cdot}$ have a coarse encoding.
This can take the form of (i) unary encoding,
(ii) entries in a fixed set of binary-encoded numbers,
or (iii) ``generalized unary encoding'' (a common generalization of
(i) and (ii)). Furthermore, to obtain
efficient algorithms, we often require that the number $d$ of
rows of $W$ is fixed.

We have four main results which below we describe and relate to the literature.
Precise statements of the theorems are given in the later sections.

\begin{enumerate}
\item[]  ({\bf Theorem \ref{thm_convexopt}}) \emph{When $\calF$ is well described, $f$ is convex (or even quasiconvex),
and $W$ has a fixed number of rows and is unary encoded or
with entries in a fixed set, we give an efficient deterministic algorithm for maximization.}
\begin{itemize}
\item[$\bullet$] As a very special case ($\calF=\{n'\times n' \mbox{ permutation
 matrices}\}$~, $W\in\Z^{d\times n'\times n'}$ unary encoded, $f$ convex), we obtain a result of
\cite{BersteinOnn2007}. The special case in which
$\calF$ is the set of characteristic vectors of common
independent sets or common bases of
a pair of matroids on $N:=\{1,2,\ldots,n\}$ appears to be new (note that in
such a case we may take $\beta=n$), even for $f$ convex.
A fortiori, also new is the special case in which
$\calF$ is the set of integer points that are in a pair of integral
polymatroids or associated base polytopes. Furthermore,
the special case in which $\calF$ is the set of characteristic vectors of
matchings of a \emph{nonbipartite} graph appears to be new, even for $f$ convex.
\end{itemize}

\item[] ({\bf Theorem \ref{MaxNorm}})  \emph{When $\calF$ is well described, $f$ is a norm, and binary-encoded $W$ is nonnegative, we give an efficient
deterministic constant-approximation algorithm for maximization.}
\begin{itemize}
\item[$\bullet$] As a very special case ($\calF=\{n'\times n' \mbox{ permutation
 matrices}\}$~, $W\in\Z^{d\times n'\times n'}$~, $f$ a $p$-norm), we obtain a result of
\cite{BersteinOnn2007}.
\end{itemize}

\item[] ({\bf Theorem \ref{approx_convex_min}}) \emph{When $\calF$ is well described, $f$ is ``ray concave''
and non-decreasing, and $W$ has a fixed number of rows
and is unary encoded or with entries in a fixed set, we give an efficient deterministic
constant-approximation algorithm for minimization.}
\begin{itemize}
\item[$\bullet$] This theorem generalizes a result of \cite{BersteinOnn2007}
which looked at the very special case of: $f\in\{$p$-norms\}$~,
$\calF:=\{n'\times n' \mbox{ permutation
 matrices}\}$~, $W\in\Z^{d\times n'\times n'}$ unary encoded.
 In addition, our more detailed analysis of $p$-norms
 further generalizes the result of \cite{BersteinOnn2007} pertaining to the
 $2$-norm.
\end{itemize}

\item[] ({\bf Theorem \ref{Random}}) \emph{When $\calF$ is the set of
characteristic vectors
of common bases of a pair
of vectorial matroids on a
common ground set, $f$ is arbitrary, and $W$ has a fixed number of rows and is unary encoded,
we give an efficient randomized algorithm for optimization.}
\begin{itemize}
\item[$\bullet$] This theorem can be contrasted with a result
in \cite{GangOfSeven2007} which established a \emph{deterministic} algorithm
for the case of a \emph{single} vectorial matroid.
\end{itemize}

\end{enumerate}

Note that all of our results extend to the generalization where
we have a binary-encoded $c\in\Z^n$ giving
a ``primary'' linear objective function $\transpose{c}x$
that we maximize, and, subject to that maximization,
we seek an optimal balancing, via
 a nonlinear function $f$~, of $d$ competing linear
 objectives $W_{i\cdot}x$~:
\[
\tag*{$P_c(\calF,\opt,f,W)$~:}
 \opt~ \left\{ ~f(Wx)~ : ~x\in \argmax\left\{\transpose{c}\tilde{x} ~:~ \tilde{x}\in \calF\right\}\right\}~.
\]
The reason is that we can find the optimal value $z^*$ of  $\transpose{c}x$ over $\calF$ using the
ellipsoid method and the separation oracle for $\calP=\conv(\calF)$~, and then
the equation $\transpose{c}x=z^*$ together with the separation oracle for
$\calP$ yields a separation oracle for the face
$\calP_c ~:=~ \calP \cap \left\{x\in\R^n ~:~ \transpose{c}x=z^*\right\}$ of $\calP$~.

Relevant background material (including polytopes, matroids,
polymatroids and linear programming) can be found in
\cite{Lee,Schrijver03,Ziegler}. Additional relevant material on
nonlinear discrete optimization can be found in
\cite{NDOBook,Onn1,OnnRothblum}.

In \S\ref{CM}, we develop an algorithm for quasiconvex combinatorial maximization.
Some of the elements that we develop for this are used in later sections as well.
In \S\ref{sec_genun}, we introduce our ``generalized unary
encoding'' scheme which we use for weight matrices.
In \S\ref{sec_fibers}, we investigate relevant properties
of the fibers (i.e., inverse images)
of points in the linear image of a polytope.
In \S\ref{sec_convexmax}, we apply properties of fibers to
give an efficient algorithm for quasiconvex combinatorial maximization.
In \S\ref{sec_approx}, we consider approximation algorithms
as a means of relaxing the assumptions that gave
us an efficient algorithm in the previous section.
We develop one algorithm for norms and another
for a very broad generalization of norms.
In \S\ref{sec_approx_convex_max}, relaxing a requirement
on the input weight matrix $W$, we
give an efficient approximation
algorithm for combinatorial maximization of norms $f$~.
In \S\ref{sec_approx_nonlinear_min}, we further exploit the
geometry of fibers to give an efficient approximation
algorithm for nonlinear combinatorial minimization, when
$W$ is nonnegative and $f$ is non-decreasing and ``ray-concave''
on the nonnegative orthant.
In \S\ref{sec_random}, for unary-encoded weight matrices $W$~,
we give an efficient algorithm for optimizing an arbitrary nonlinear
function $f$~, over the set $\calF$ of characteristic vectors of
common bases of a pair of vectorial matroids on a common ground set.

\section{Quasiconvex combinatorial maximization}\label{CM}

\subsection{Generalized unary encoding of weights}\label{sec_genun}
General binary encoding of weights is typically too parsimonious to
allow for the construction of theoretically efficient algorithms for nonlinear
discrete optimization.
So we turn our attention to less parsimonious encodings.
We  consider weights $W_{i,j}$ of the form
\[
     W_{i,j} = \sum_{k=1}^p \delta^k_{i,j}  a_k~  ,
\]
 with integer $p\ge 1$ fixed, the distinct positive integers $a_k$ being binary encoded, but the
integers $\delta^k_{i,j}$ (unrestricted in sign) being unary encoded.
It is convenient to think of unary-encoded
matrices $\delta^k=((\delta^k_{i,j}))\in\Z^{d\times n}$~, for $1\le k\le p$~,
and binary encoded $a\in\Z^p$~, with $a_k>0$~, for $1\le k\le p$~. Then we
have $W:=((W_{i,j}))=\sum_{k=1}^p a_k \delta^k$~.

We have the following special cases:
\begin{enumerate}
\item\label{unary_} {\bf Unary-encoded weights:} With $p=1$ and $a_1=1$~, we get the ordinary model of
unary-encoded $W=\delta^1$~.
\item\label{a1ap_} {\bf Binary-encoded \mbox{\boldmath$\{0,a_1 , a_2 ,..., a_p\}$}-valued weights:}
With $\delta^k\ge 0$ for
all $k$~, and $\sum_{k=1}^p \delta^k_{i,j} \le 1$~, for
all $i,j$~, we get the case of all $W_{i,j}$ in the
set $\{0,a_1 , a_2 ,..., a_p\}$ having binary-encoded elements~.
\item\label{binary_} {\bf 0/1-valued weights:} This is, of course, the important common
special case of \ref{unary_} and \ref{a1ap_}.
\end{enumerate}
Because \ref{unary_} is a much more common generalization of \ref{binary_} than is \ref{a1ap_},
we refer to our general setting as \emph{generalized unary encoding}.
We do note that there are cases where, with respect to
a particular nonlinear combinatorial optimization problem,
one algorithm is efficient for case \ref{a1ap_}
but not for case \ref{unary_}, and for another algorithm vice versa
(see  \cite{GangOfSeven2007,LOWFrobenius} for example).
Furthermore, there are broad cases where
with binary encoding a problem is provably intractable,
but we will see that we can
construct an efficient algorithm when the
weights have a generalized unary encoding.

\subsection{Fibers}\label{sec_fibers}

We begin with some definitions.
For $W\in \Z^{d\times n}$ and polytope $\calP\subset \R^n$~,
we let $W{\calP}$ denote the image of
$\calP$ under the linear map $x\mapsto Wx$~.
Clearly $W{\calP}\subset \R^d$
is also a polytope. For $u\in\R^d$~,
we define
\[
W^{-1}u := \left\{x\in \R^n ~:~ Wx=u \right\}~.
\]
For $u\in\R^d$~, we refer to $(W^{-1}u) \cap \calP$
as the $W_{\calP}$-fiber of $u$~.
If $W$ and
$\calP$ are clear from context, then
we may just say ``the fiber of $u$~.''
For $u\in \R^d$~, we note that $u\in W\calP$
if and only if the $W_{\calP}$-fiber of $u$ is nonempty.

First, we point out how we can optimize
efficiently a linear function on the
fiber of any $u\in\R^d$~.
In particular, we can also determine nonemptiness
of such a fiber.

\begin{lemma}\label{fiber_LP}
Assume that we have
a separation oracle for the polytope $\calP\subset \R^n$~.
Let $c\in\Z^n$~, $W\in\Z^{d\times n}$ and  $u\in\R^d$
be binary encoded.
Then we can solve the linear-objective optimization
problem
\[
\max
\left\{
\transpose{c}x\ : x\in (W^{-1}u) \cap {\calP}
 \right\}
\]
in polynomial time. In particular,
we can test nonemptiness of the $W_{\calP}$-fiber
of $u$
in polynomial time.
\end{lemma}

\begin{proof}
The problem is just the
linear program
\[
\max
\left\{
\transpose{c}x\ : Wx=u~,\ x\in\calP
 \right\}~,
\]
which we can solve in polynomial time via
the ellipsoid method, as we assume that we have available
a separation oracle for $\calP$~.
\qed
\end{proof}

Unfortunately, even though $\calP$ is assumed to
have integer vertices (it is the convex hull of
the finite set $\calF$ of integer points), the $W_{\calP}$-fiber
of $u$ may have some integer vertices and
some fractional ones (after all,
we are just intersecting $\calP$ with a linear
subspace). Therefore, an optimal extreme
point solution of the linear program of Lemma \ref{fiber_LP}
may be at a fractional vertex of the $W_{\calP}$-fiber
of $u$~. In short, fibers of arbitrary
$u\in W\calP$ are not well behaved. However,
we will see that fibers of $u\in \vert(W\calP)$
are better behaved. Specifically, we will now demonstrate that
if $u\in \vert(W\calP)$~, then we can find
efficiently an integer optimizer of
any linear  function on the
$W_{\calP}$-fiber of $u$~.

\begin{lemma}\label{fiber_IP}
Assume that we have
a separation oracle for the polytope $\calP\subset \R^n$~.
Let $c\in\Z^n$~, $W\in\Z^{d\times n}$ and $u\in \vert(W\calP)$
be binary encoded.
Then we can solve the optimization
problem
\[
\max
\left\{
\transpose{c}x\ : x\in (W^{-1}u) \cap \calP \cap \Z^n
 \right\}
\]
in polynomial time.
\end{lemma}

\begin{proof}
It is well known that (i) the image $\calQ$ of a polytope $\calP$
under a linear (even affine) map is a polytope, and (ii) the
preimage (in $\calP$) of every face of $\calQ$ is a face of
$\calP$~. In particular, the
preimage (in $\calP$) of every vertex of $\calQ$ is a face of
$\calP$~.

Therefore, if $u\in \vert(W\calP)$~, then the $W_{\calP}$-fiber of $u$
is a face of $\calP$~, and so such a fiber
has itself integer vertices (because $\calP$ does). Therefore,
an extreme-point solution of the linear program
\[
\max
\left\{
\transpose{c}x\ : Wx=u~,\ x\in\calP
 \right\}
\]
will be integer.
\qed
\end{proof}

The following result is not needed for our purposes,
but it is interesting in its own right.

\begin{proposition}\label{fiber_descrip}
Assume that we have
$\calP=\{x\in\R^n : Ax\le b\}$~, for some
binary encoded $A\in \Z^{m\times n}$ and $b\in \Z^m$~.
Let $W\in\Z^{d\times n}$ and $u\in \vert(W\calP)$
be binary encoded,
Then we
can compute an inequality description
of the $W_{\calP}$-fiber of $u$ as a face
of $\calP$ in polynomial time.
\end{proposition}

\begin{proof}
For each row $1\le i\le m$~,
we simply solve the linear program
\[
z_i:=\min_{x\in\R^n}
\left\{
A_{i\cdot} x\ : Wx=u~,\ Ax\le b
 \right\}~.
\]
Then
\[
(W^{-1}u ) \cap \calP = \left\{x\in \R^n\ :\ A_{i\cdot}x
\left\{
\hskip-3pt
\begin{array}{c}
\le\\[-2pt]
=
\end{array}
\hskip-3pt
\right\}
b_i
\mbox{ if } z_i
\left\{
\hskip-3pt
\begin{array}{c}
<\\[-2pt]
=
\end{array}
\hskip-3pt
\right\}
b_i~, \mbox{ for } 1\le i\le m
\right\}~.
\]
\qed
\end{proof}

Next, we demonstrate that under an appropriate encoding
of the data, the set of vertices of $W\calP$
can be computed efficiently.

\begin{lemma}\label{image_vert}
Let $W\in\Z^{d\times n}$
be given.
We assume that $\calF$ is a finite subset of $Z^n_+$~, and further that
 $\calF\subset\calS(n,\beta):=\{ x\in \R_+^n\ :\ \transpose{\vecone}x \le \beta\}~,$
 for some $\beta \in \Z_+$~. Then, for fixed $d$~, we can compute the
 vertices of $W\calP$ in time polynomial in $n$~, $\beta$~, and
 the length of the generalized unary encoding of $W$~.
\end{lemma}

\begin{proof}

For $x\in\calF$~,
consider $u:=Wx$~. For $1\le i\le d$~, we have
\[
\begin{array}{rrl}
u_i ~:=~ W_{i \cdot} x
~=~ \sum_{j=1}^n W_{i,j} x_j
~=~ \sum_{j=1}^n \left(\sum_{k=1}^p \delta^k_{i,j}  a_k\right) x_j
~=~ \sum_{k=1}^p \left( \sum_{j=1}^n \delta^k_{i,j} x_j \right) a_k~.
\end{array}
\]

Observe that the nonnegative integer coefficients $\sum_{j=1}^n \delta^k_{i,j} x_j$
of the $a_k$ are bounded by a polynomial in the unary $\beta$ and the unary $\delta^k_{i,j}$~:
Letting $\omega:= \max \left|\delta^k_{i,j}\right|$~, we observe that
\[
{\textstyle
\left|
\sum_{j=1}^n \delta^k_{i,j} x_j
\right|
\le
\left\|  \delta^k_{i,j} \right\|_{\infty}
\left\| x \right\|_{1}
\le
\omega\beta ~.}
\]
This means that there are only a polynomial number of possible
values for each $u_i$ ($p$ is fixed), and then (with $d$ fixed)
only a polynomial number of images $u=Wx$ of $x\in\calF$~.

Concretely, for $k=1,\ldots,p$~, we let $\hat{u}^k$ range over
\[
\{0,\pm 1,\ldots,\pm\omega\beta\}^d~\subset\Z^d~.
\]
For each of the $(2\omega\beta+1)^{pd}$ choices of
$(\hat{u}^1,\ldots,\hat{u}^p)$~, we let
$u:=\sum_{k=1}^p a_k \hat{u}^k$~. For each
such $u$ we check if $u\in W\calP$  using the algorithm of Lemma \ref{fiber_LP},
and so we obtain a set
$U\subset \Z^d$~, comprising these points $u$ having a nonempty fiber,
having a polynomial number of elements.

Now, because $U$ is contained in $W\calP$ and
contains $\vert(W\calP)$~, it is
clear that $W\calP=\conv(U)$~.

Finally, because the dimension $d$ is fixed,
we can compute the vertices of $\conv(U)$
in polynomial time.
\qed
\end{proof}

\subsection{Exact algorithm for quasiconvex maximization}\label{sec_convexmax}

In this subsection, we take the facts that we have already built up to
give an efficient algorithm to maximize a function that
is  quasiconvex
(that is, its ``lower level sets'' $\{z\in\R^d ~:~ f(z)
\le \tilde{f}\}$ are convex subsets of $\R^d$~, for all $\tilde{f}\in
\R$~; see \cite{ADSZ}, for example). Of course, the algorithm
also applies to ordinary convex functions.

\begin{theorem} \label{thm_convexopt}
Suppose that we are given
$W\in \Z^{d\times n}$~, quasiconvex function $f:\R^d\rightarrow\R$
specified by a comparison oracle, and a finite subset $\calF$
of $\Z^n_+$~. We further assume that
 $\calF\subset\calS(n,\beta):=\{ x\in \R_+^n\ :\ \transpose{\vecone}x \le \beta\}~,$
 for some $\beta \in \Z_+$~, and that
 we have a separation oracle for $\calP:=\conv(\calF)$~.
 Then, if $d$ is fixed, we can solve
 $P(\calF,\max,f,W)$ in time polynomial in $n$~, $\beta$~,
 and the
 size of the generalized unary encoding of $W$~.
\end{theorem}

\begin{proof}
Because $f(u)$ is quasiconvex on the polytope $W\calP$
and $f_W(x):=f(Wx)$ is
quasiconvex on the polytope $\calP$~, and the maximum of a quasiconvex function on
a polytope is attained at a vertex of the polytope, we have
\[
\begin{array}{l}
\max \left\{ f(Wx)\ :\ x\in\vert(\calP)\right\}
~~=~~ \max \left\{ f(Wx)\ :\ x\in\calP\right\}\\
\quad = \max \left\{ f(u)\ :\ u\in W\calP\right\}
~~=~~ \max \left\{ f(u)\ :\ u\in\vert(W\calP)\right\}~.
\end{array}
\]

First, via Lemma \ref{image_vert}, we compute efficiently $\vert(W\calP)$~.
Next, using the comparison oracle for $f$~, we can compare values of
$f(u)$~, for all $u\in\vert(W\calP)$~,
to obtain the best $u$~.
Then, for the best $u\in\vert(W\calP)$~,
via Lemma \ref{fiber_LP}, using an arbitrary linear objective
vector $c\in\Z^n$ in general position,
we efficiently find an extreme point $x^u$ of the $W_{\calP}$-fiber of $u$~.
By Lemma \ref{fiber_IP}, such a point $x^u$ will be integer and will thus solve
$P(\calF,\max,f,W)$~.
\qed
\end{proof}

\section{Approximation algorithms for norms and related functions}\label{sec_approx}

In this section, we consider approximation algorithms
as a means of relaxing the assumptions that gave
us an efficient algorithm in the previous section.
In \S\ref{sec_approx_convex_max}, we provide an efficient approximation algorithm for
combinatorial norm maximization.
In \S\ref{sec_approx_nonlinear_min}, we provide an efficient approximation algorithm
for combinatorial minimization. This latter algorithm applies to
a very broad generalization of norms.

\subsection{Approximation algorithm for norm maximization}\label{sec_approx_convex_max}

The following theorem provides an approximation algorithm for
maximizing a norm $f$~, that runs in time that is polynomial
even in the bit size of the weights $W^k_{i,j}$ and even if $d$ is
variable. The approximation ratio depends on the (positive) constants
for the standard equivalence of the norm $f$ with the infinity norm:
\[
C_f ~ \| u \|_\infty ~ \le  ~ f(u) ~  \le ~ C^f ~ \| u \|_\infty~.
\]

Note that appropriate constants $C_f$ and $C^f$ can be concretely
described, and in some cases efficiently calculated.
Letting $\{e_1,\ldots,e_d\}$ denote the standard basis of $\R^d$~,
a valid choice of $C^f$ is the sum of $f$ evaluated on each of the basis
elements (i.e., $C_f:=\sum_{i=1}^d f(e_i)$).
Furthermore,
the best $C_f$ is the minimum of the convex function $f$ on the unit sphere in the
$\infty$-norm (i.e., $C_f:= \min \{f(u) ~:~ \|u\|_\infty = 1\}$~.
This is not a convex optimization problem because we are restricted to the
\emph{boundary} of the (unit) ball (of the $\infty$-norm). If however we
take $d$ to be constant or even say $O(\log~n)$~, then
we can solve $2^d$ convex minimization problems (via say an ellipsoid method),
one for reach facet of the unit ball of the $\infty$-norm,
in order to efficiently calculate a valid (in fact, best) value of
$C_f$~.
Finally, for the special case of $f(u):=\|u\|_p$~,
our algorithm finds a
$d^{1/p}$-approximate
solution.

\begin{theorem}\label{MaxNorm}
Suppose that we are given nonnegative $W\in \Z^{d\times n}$~, and $\calF$
a finite subset of $\Z^n_+$~. We further assume that
 $\calF\subset\calS(n,\beta):=\{ x\in \R_+^n\ :\ \transpose{\vecone}x \le \beta\}$~,
 for some $\beta \in \Z_+$~, and that
 we have a separation oracle for $\calP:=\conv(\calF)$~.
Then, for any norm $f$ given by a comparison oracle,
there is an algorithm that determines a $\left(C^f/C_f\right)$-approximate
solution to $P(\calF,\max,f,W)$~, in time that is polynomial in $d$~,$n$~, $\beta$~,
and $\max\left\lceil \log w^k_{i,j}\right\rceil$~.
Moreover, for $f(u):=\|u\|_p$~,
for any $p$ satisfying $1\leq p\leq\infty$~,
our algorithm determines a
$d^{1/p}$-approximate
solution.
\end{theorem}

\begin{proof}
First, let $f$ be an arbitrary norm.
The algorithm
is the following: For $i=1,\dots,d$~, solve the linear-programming
problem
\[
\max\{W_{i\cdot} x ~:~ x\in \calP \}~ ,
\]
obtaining an optimal vertex $x^i$ of $\calP$~, and let $u^i$ be its $W$-image.
Then output
$x^r$ such that $\|u^r\|_p=\max_{i=1}^d\|u^i\|_p$~.

We now show that this provides the claimed approximation. Let $s$
satisfy $\|u^s\|_\infty=\max_{i=1}^d\|u^i\|_\infty$~. First, we
claim that any $\tilde{u}\in\vert(W\calP)$ satisfies $\|\tilde{u}\|_\infty\leq
\|u^s\|_\infty$~. To see this, recall that
$W$ and hence the $u^i$ are nonnegative; choose any point $\tilde{x}\in
(W^{-1}\tilde{u}) \cap \calF$ (hence $\tilde{u}=W\tilde{x}$), and let $t$
satisfy $\tilde{u}_t=\|\tilde{u}\|_\infty=\max_{i=1}^d \tilde{u}_i$~.  Then, as claimed,
we get
\[
\|\tilde{u}\|_\infty \ \ = \ \ \tilde{u}_t \ \ =\ \ W_{t\cdot} \tilde{x}
\ \ \leq \ \ \max\{W_{t\cdot} x\,:\, x\in\calF\} \ \ = \ \ W_{t\cdot} x^t \ \ =
\ \ u^t_t \ \ \leq \ \ \|u^t\|_\infty \ \ \leq\ \ \|u^s\|_\infty \
\ .
\]

Let $x^*$ be an optimal solution, and
let $u^*$ be its $W$-image. As norms are (quasi)convex,
we can without loss of generality take $x^*$ to be a vertex
of $\calP$ and $u^*$ to be a vertex of $W\calP$~.

Then we have the following inequalities:
\begin{eqnarray*}
f(Wx^*)
& = &
f(u^*)
\quad \leq \quad
C^f ~ \|u^*\|_\infty
\quad \leq \quad
C^f ~ \|u^s\|_\infty \\
& \leq &
\left(C^f/C_f \right) ~ f(u^s)
\quad \leq \quad
\left(C^f/C_f \right)  ~ f(u^r)
\quad = \quad
\left(C^f/C_f \right)  f(Wx^r)~.
\end{eqnarray*}

Finally, observing that
for $p$ satisfying $1\leq p\leq\infty$~,
we have  $C_f=1$ and $C^f=d^{1/p}$~,
we conclude that for $p$-norms
we get a $d^{1/p}$-approximate
solution.
\qed
\end{proof}

\subsection{Approximation algorithm for ray-concave minimization}\label{sec_approx_nonlinear_min}

In this subsection, using some of the methodology that we already developed, we obtain
an efficient approximation algorithm for nonlinear combinatorial minimization. Our algorithm
applies to certain functions that generalize norms via, what will at first seem
paradoxical, a type of concavity.

A function $f:\R^d_+\rightarrow \R$ is
\emph{ray-concave} if
\begin{enumerate}
\item[(i)] $\lambda f(u) \le f(\lambda u) \mbox{ for } u\in\R^d_+$~,\ $0\le \lambda\le 1$~;
\end{enumerate}
Analogously $f:\R^d_+\rightarrow \R$ is
\emph{ray-convex} if
\begin{enumerate}
\item[(i')] $\lambda f(u) \ge f(\lambda u) \mbox{ for } u\in\R^d_+$~,\ $0\le \lambda\le 1$~;
\end{enumerate}
Naturally, we say that $f:\R^d_+\rightarrow \R$ is \emph{homogeneous} if it is both ray-concave and ray-convex.

A function $f:\R^d_+\rightarrow \R$~ is \emph{non-decreasing} if
\begin{enumerate}
\item[(ii)] $f(u)\le f(\tilde{u})~, \mbox{ for } u,\tilde{u}\in\R^d_+~,\ u\le\tilde{u}$~.
\end{enumerate}

We are mainly interested in ray-concave non-decreasing functions $f:\R^d_+\rightarrow\R$~.
In particular, by $(i)$ with $\lambda=0$~, we have $f(\veczero)\ge 0$~,
and hence by $(ii)$~, we have $f(u)\ge 0$ for all $u\in\R^d_+$~.

Notice that ordinary concavity of a function $f$ has the special case:
\[
\lambda f(u) + (1-\lambda) f(\veczero) \le f(\lambda u + (1-\lambda)\veczero)~, \mbox{ for }
u\in\R^d_+~,\ 0\le \lambda\le 1~,
\]
and so if $f$ is concave with $f(\veczero)=0$~, then it is ray-concave
--- thus justifying our terminology.

Notice further that if $f$ is a norm on $\R^d$~, then it of course non-decreasing on
$\R^d_+$~; additionally, it is not only ray-concave, but it is more strongly homogeneous
(and more of course). In general, homogeneity already implies $f(\veczero)=0$~.
A concrete example of a ray-convex function that is
not generally homogeneous is $h(u):=\prod_{i=1}^d u_i$~.

There are other interesting and useful functions that are ray-concave.
Suppose that ray-concave $g: \R^d_+ \rightarrow \R$
and ray-convex $h: \R^d_+ \rightarrow \R$
satisfy
\begin{enumerate}
\item[(iii)] $g(u+v)-g(u)\ge h(u+v)-h(u)$~, for all $u,v\in \R^d_+$~.
\end{enumerate}
That is, $g$ grows no slower than $h$ on $\R^d_+$~.
Then $f(u):=g(u)-h(u)$ is of course non-decreasing.
Moreover, if $g$ is ray-concave on $\R^d_+$ and $h$ is ray-convex on $\R^d_+$~,
then $f$ is also ray-concave on $\R^d_+$~.

As a concrete and natural example, consider $g(u):=\| u \|_1$ and $h(u):=\| u \|_p$
for any integer $p\ge 1$ or infinity.
Then $f(u):=g(u)-h(u)=\| u \|_1 - \| u \|_p$ is ray-concave and non-decreasing
on $\R^d_+$~.
Notice that already for $d=2$ and $p=\infty$~, $f(u)$ is not a norm --- indeed,
for this case $f(u)=\min(u_1,u_2)$~.

\begin{theorem}\label{approx_convex_min}
Suppose that nonnegative $W\in \Z^{d\times n}$ is given, $d$ is fixed,
$f:\R^d_+\rightarrow\R$ is a ray-concave non-decreasing function
given by a comparison oracle, and $\calF$ is
a finite subset of $\Z^n_+$~. We further assume that
 $\calF\subset\calS(n,\beta):=\{ x\in \R_+^n\ :\ \transpose{\vecone}x \le \beta\}~,$
 for some $\beta \in \Z_+$~, and that
 we have a separation oracle for $\calP:=\conv(\calF)$~.
 Then we can compute a $d$-approximate
 solution to $P(\calF,\min,f,W)$ in time
 polynomial in $n$~, $\beta$~, and the size of the generalized unary encoding of $W$~.
Moreover, for $f(u):=\|u\|_p$ (the $p$-norm), $1\le p\le \infty$~,
the algorithm actually determines a  $d^{1/q}$-approximate solution, where
$1/p+1/q=1$ (with obvious interpretations when a denominator is $\infty$)~.
\end{theorem}

\begin{proof}
Because we are minimizing a function that need not be concave, it may well
be the case
that the optimal solution of $P(\calF,\min,f,W)$ is \emph{not}
realized at an extreme point of $\calP$ and consequently, it may be that
$\min\{ f(u) \ :\ u\in  W\calP\}$ is
\emph{not} realized at an extreme point of $W\calP$~.
Nonetheless, we will focus on extreme points of
$W\calP$~.

Apply the algorithm of Lemma \ref{image_vert} so as to construct the
set of vertices of $W\calP$~. Using the comparison
oracle of $f$~, identify a vertex $\hat{u}\in
\vert(W\calP)$ attaining minimum value $f(u)$~. Now apply the
algorithm of Lemma~\ref{fiber_LP} to $\hat{u}$ and, as guaranteed by
Lemma \ref{fiber_IP}, obtain an integer $\hat{x}$ in
the fiber of $\hat{u}$~, so that $\hat{u}=W\hat{x}$~. Output
the point $\hat{x}$~.

We now show that this provides the claimed approximation. Let $x^*$
be an optimal integer point, and let $u^*:=W x^*$ be its image.
Let $u'$ be a
point on the boundary of $W\calP$ satisfying $u'\leq u^*$~.
By
Carath\'eodory's theorem (on a facet of
$W\calP$ containing $u'$),\
$u'$ is a convex
combination $u'=\sum_{i=1}^r\lambda_i u^i$ of some $r\leq d$ vertices of
$W\calP$ for some coefficients $\lambda_i\geq 0$ with $\sum_{i=1}^r\lambda_i =1$~.
Let $t$ be an index for which $\lambda_t=\max_i \{ \lambda_i\}$~.
Then $\lambda_t\geq{1\over r}\sum_{i=1}^r\lambda_i=1/r\geq 1/d$~.

Because $W$ is nonnegative, we find that so are $u'$
and the $u^i$~, and hence we obtain
\begin{align*}
 & f(W\hat{x}) \quad = \quad
f(\hat{u}) \quad \leq \quad
f(u^t) \quad \leq \quad
d\lambda_t\cdot f(u^t) \quad \leq \qquad
d\cdot f(\lambda_t u^t) \quad\\
 &\qquad \leq \quad {\textstyle d\cdot f(\sum_{i=1}^r\lambda_i u^i) } \quad = \quad
d\cdot f(u') \quad \leq \quad
d\cdot f(u^*) \quad = \quad
d\cdot f(Wx^*)~.
\end{align*}
This proves that $\hat{x}$ provides a $d$-approximate solution.

Now consider the case of the $p$-norm
$f(u)=\|u\|_p$~. First, we note that the first part already covers the
case of $p=\infty$~. So we confine our attention to $p<\infty$~.
We complete now the proof for $1<p<\infty$~, and then it is an easy observation that
the proof still goes through (in a simplified form) establishing also
that for $p=1$ we get a $1$-approximate (i.e., optimal) solution.

By the H\"older inequality
($|\transpose{a} b| \le \|a\|_p \|b\|_q$~, for $1/p+1/q=1$),
we have
\[
\textstyle
1
~=~
\sum_{i=1}^r 1\cdot\lambda_i
\leq
\left(\sum_{i=1}^r 1^q \right)^{1/q} \left(\sum_{i=1}^r\lambda_i^p\right)^{1/p}
~=~
r^{1/q} \left(\sum_{i=1}^r\lambda_i^p\right)^{1/p}
~\leq~
d^{1/q} \left(\sum_{i=1}^r\lambda_i^p\right)^{1/p}~.
\]
Find $s$ with $\|u^s\|_p=\min_i \{\|u^i\|_p\}$ and recall that the $u^i$
are nonnegative. We then have
\begin{align*}
\textstyle
& f^p(W\hat{x})
\quad =  \quad
\|\hat{u}\|_p^p
\quad\leq\quad
\|u^s\|_p^p
\quad\leq\quad
\textstyle d^{p/q} \left(\sum_{i=1}^r\lambda_i^p\right) \|u^s\|_p^p
\quad\leq\quad
\textstyle d^{p/q} \sum_{i=1}^r \left( \lambda_i^p \|u^i\|_p^p \right) \\
& \qquad \leq
\textstyle d^{p/q} \left\|\sum_{i=1}^r\lambda_i u^i\right\|_p^p
\quad = \quad
d^{p/q} \|u'\|_p^p
\quad \leq\quad
d^{p/q}\|u^*\|_p^p
\quad =\quad
d^{p/q} f^p(Wx^*)~,
\end{align*}
which proves that in this case, as claimed, $\hat{x}$ provides
moreover a $d^{1/q}$-approximate solution.
\qed
\end{proof}

\section{Randomized algorithm for vectorial matroid intersection}\label{sec_random}

The algorithm of \S\ref{sec_convexmax} can be used to solve efficiently
the problem $P(\calF,\max,f,W)$~, when $\calF$ is the set of characteristic
vectors of sets that are common independent sets or
bases of a pair of matroids on a common ground set,
and $f$ is convex.
Such an algorithm only needs an independence oracle for each of the matroids.
Even better, we can do likewise when $\calF$ is the set of integer points that
are in a pair of integral polymatroids or associated base polytopes. For polymatroids,
we only need an oracle that evaluates the associated submodular functions.

In this section, we relax the assumption of $f$ being convex.
Indeed, we give a \emph{randomized} algorithm for arbitrary $f$
presented by a comparison oracle. By this we mean an algorithm that
has access to a random bit generator, and on any input, it outputs the
optimal solution with probability at least half. We restrict our
attention to unary encoded $W$ and $\calF$ being the set of
characteristic vectors of common bases of a pair of vectorial
matroids $M_1$ and $M_2$ on the common ground set
$\{1,2,\ldots,n\}$~. Our broad generalization of the standard
linear-objective matroid intersection problem is defined as follows.

\vskip.2cm\noindent {\bf Nonlinear Matroid Intersection}. Given two
matroids $M_1$ and $M_2$ on the common ground set $N=\{1,\ldots,n\}$~,
integer weight matrix $W\in \Z^{d\times n}$ and an arbitrary
function $f:\R^d\rightarrow\R$ specified by a comparison oracle~.
Let $\calF$ be the set of characteristic vectors of common bases of $M_1$ and $M_2$~.
Find an $x\in{\cal F}$ maximizing (or minimizing)
$f(Wx)$~.
\vskip.2cm

When $M_1$ and $M_2$ are vectorial matroids, they can be represented
by a pair of $r \times n$ matrices. We assume, without loss of
generality, that $r$ is the common rank of $M_1$ and $M_2$. Note
that in this case $\calF$ satisfies $\calF\subseteq I_{n,r}$~, where
$I_{n,r}$ is the set of vectors in $\R^n$ with exactly $r$ nonzero
entries which are equal to 1. Thus our problem consists on finding
an $x$ that maximizes (minimizes) $f(Wx)$~, where $x$ is the
characteristic vector corresponding to the set of column indices
defining a common base. We assume the matroids have at least one
common base which can be efficiently checked by the standard matroid
intersection algorithm.

By adding to each $W_{i,j}$ a suitable positive integer $v$ and
replacing the function $f$ by the function that maps each $u\in\R^d$
to $f(u_1-rv,\dots,u_d-rv)$ if necessary, we may and will assume
without loss of generality throughout this section that the given
weight matrix is nonnegative, $W\in \Z_+^{d\times n}$~.

In this section we will be working with polynomials with integer
coefficients in the $n+d$ variables $a_{j}$~, $j=1,\dots,n$ and
$b_i$~, $i=1,\dots,d$~. Define a vector $\gamma$ in $\R^n$ whose
entries are monomials by
\vskip-5pt
\[
\gamma_j ~:=~ a_j\prod_{i=1}^db_i^{W_{i,j}},~ j=1,\ldots,n~.
\]
\vskip-5pt
For each vectors $x\in \Z_+^{n}$ and $u\in \Z_+^{d}$~, the
corresponding monomials are
\vskip-5pt
\[
a^x:=\prod_{j=1}^n a_{j}^{x_{j}} \quad,\qquad b^u:=\prod_{i=1}^d b_i^{u_i}~.
\]
\vskip-5pt Let $M_1^x$ and $M_2^x$ be the $r\times r$ matrices
comprising the $r$ columns whose indices correspond to the nonzero
entries in $x$~. Finally, for each $u\in \Z_+^{d}$ define the
following polynomial in the variables $a=(a_{j})$ only.
\[
g_u(a)\quad:=\quad\sum
\left\{\det(M_1^x)\cdot\det(M_2^x)a^x\,:\,x\in \calF~,\ \ W
x=u\right\}~.
\]
 Define $U:=W\calF=\{Wx: x\in \calF\}$ as the set of
points that are the projection of $\calF$ under the weight matrix
$W$~. Let $M_1(\gamma)$ be a matrix whose $j$-th column is $M_1^j
\gamma_j$~. We then have the following identity in terms of the
$g_u(a)$~.
\begin{align*}
&\det(M_1(\gamma)\transpose{M}_2)
~=~ \sum_{x\in I_{n,r}}\det(M_1^x(\gamma))\cdot\det(M_2^x)
~=~\sum_{x\in I_{n,r}}\det(M_1^x)\cdot\det(M_2^x)\prod_{j=1}^n \gamma_j^{x_j} \displaybreak[0]\\
&\quad =\sum_{x\in \calF}\det(M_1^x)\cdot\det(M_2^x)\prod_{j=1}^n
\gamma_j^{x_j}~=~\sum_{x\in
\calF}\det(M_1^x)\cdot\det(M_2^x)\prod_{j=1}^n(a_j\prod_{i=1}^db_i^{W_{i,j}})^{x_j}\displaybreak[0]\\
&\quad=\sum_{x\in \calF}\det(M_1^x)\cdot\det(M_2^x)\prod_{j=1}^n
a_j^{x_j}\prod_{j=1}^n\prod_{i=1}^d b_i^{W_{i,j} x_j}~=~\sum_{x\in
\calF}\det(M_1^x)\cdot\det(M_2^x)\prod_{j=1}^n
a_j^{x_j}\prod_{i=1}^d b_i^{W_i\cdot x} \displaybreak[0]\\
&\quad= \sum_{x\in \calF}\det(M_1^x)\cdot\det(M_2^x)a^xb^{W\cdot
x}~=~\sum_{u\in U}\sum_{{x\in \calF~: \atop W x=u}}
\det(M_1^x)\cdot\det(M_2^x)a^xb^{W\cdot x} ~=~ \sum_{u\in U}
g_u(a)b^{u} ~.
\end{align*}

 Next we consider integer
substitutions for the variables $a_j$~. Under such substitutions, each
$g_u(a)$ becomes an integer, and $\det(M_1(\gamma)\transpose{M}_2)=\sum_{u\in
U}g_u(a)b^u$ becomes a polynomial in the variables $b=(b_i)$ only.
Given such a substitution, let ${U(a)}:=\{u\in U\,:\,g_u(a)\neq 0\}$
be the {\em support} of $\det(M_1(\gamma)\transpose{M}_2)$~, that is, the set
of exponents of monomial $b^u$ appearing with nonzero coefficient in
$\det(M_1(\gamma)\transpose{M}_2)$~.

The next proposition concerns substitutions of independent identical
random variables uniformly distributed on
$\{1,2,\dots,s\}$~, under which $U(a)$ becomes a random subset $\widehat
U\subseteq U$~.

\begin{proposition}\label{Probability}
Suppose that independent uniformly-distributed
random variables  on
$\{1,2,\dots,s\}$ are substituted for the $a_j$~, and let ${\widehat
U}:=\{u\in U\,:\,g_u(a)\neq 0\}$ be the random support of
$\det(M_1(\gamma)\transpose{M}_2)$~. Then, for every $u\in U=\{W
x\,:\,x\in\calF\}$~, the probability that $u\notin\widehat U$ is at most
${r/s}$~.
\end{proposition}

\begin{proof}
Consider any
$u\in U$~, and consider $g_u(a)$ as a polynomial in the variables
$a=(a_j)$~. Because $u=Wx$ for some $x\in\calF$~, there is as least one
term $\det(M_1^x)\cdot\det(M_2^x)a^x$ in $g_u(a)$~. Because distinct
$x\in\calF$ give distinct monomials $a^x$~, no cancelations
occur among the terms $\det(M_1^x)\cdot\det(M_2^x)a^x$ in $g_u(a)$~.
Thus, $g_u(a)$ is a nonzero polynomial of degree $r$~. The claim now
follows from a lemma of Schwartz \cite{Sch} stating that the
substitution of independent uniformly-distributed random variable
on $\{1,2,\dots,s\}$ into a nonzero multivariate
polynomial of degree $r$ is zero with probability at most $r/s$~.
\qed
\end{proof}

The next lemma establishes that, given $a_j$~, the support $U(a)$ of
$\det(M_1(\gamma)\transpose{M}_2)$ is polynomial-time computable.

\begin{lemma}\label{Interpolation}
For every fixed $d$~, there is an algorithm that,
given $W\in\Z_+^{d\times n}$~, and substitutions
$a_{j}\in\{1,2,\dots,s\}$~, computes $U(a)=\{u\in U\,:\,g_u(a)\neq
0\}$ in time polynomial in $n$~, $\max W_{i,j}$ and $\lceil \log s
\rceil$~.
\end{lemma}

\begin{proof}
The proof is based on interpolation.
For each $u$~, let $g_u:=g_u(a)$ be the fixed integer obtained by
substituting the given integers $a_{j}$~. Let $z := r \max W_{i,j}$
and $Z:=\{0,1,\dots, z\}^d$~. Note that $z$ and $|Z|=(1 + z)^d$ are
polynomial in $n$ and the unary encoding of $W$~. Then
$U(a)\subseteq U\subseteq Z$~, and hence
$\det(M_1(\gamma)\transpose{M}_2)=\sum_{u\in Z}g_ub^u$ is a
polynomial in $d$ variables $b=(b_i)$ involving at most
$|Z|=(z+1)^d$ monomials. For $t=1,2,\dots,(z+1)^d$~, consider the
substitution $b_i:=t^{(z+1)^{i-1}}$~, $i=1,\dots,d$~, of suitable
points on the moment curve. Under this substitution of $b$ and a
given substitution of $a$~, the matrix
$M_1(\gamma(t))\transpose{M}_2$ becomes an integer matrix and so
its determinant $\det(M_1(\gamma(t))\transpose{M}_2)$ becomes an
integer number that can be computed in polynomial time by Gaussian
elimination. So we obtain the following system of $(z+1)^d$
equations in $(z+1)^d$ variables $g_u$~, $u\in Z=\{0,1,\dots,
z\}^d$~,
\[
\det(M_1(\gamma(t))\transpose{M}_2)\quad=\quad\sum_{u\in
Z}g_u\prod_{k=1}^d b_k^{u_k} \quad=\quad\sum_{u\in Z}t^{\sum_{k=1}^d
u_k (z+1)^{k-1}}\cdot g_u \,,\quad\quad t=1,2,\dots,(z+1)^d~.
\]
As $u=(u_1,\dots,u_d)$ runs through $Z$~, the sum $\sum_{k=1}^d u_k
(z+1)^{k-1}$ attains precisely all $|Z|=(z+1)^d$ distinct values
$0,1,\dots,(z+1)^d-1$~. This implies that, under the total order of
the points $u$ in $Z$ by increasing value of $\sum_{k=1}^d u_k
(z+1)^{k-1}$~, the vector of coefficients of the $g_u$ in the
equation corresponding to $t$ is precisely the point
$(t^0,t^1,\dots,t^{(z+1)^d-1})$ on the moment curve in
$\R^Z\simeq\R^{(z+1)^d}$~. Therefore, the equations are linearly
independent, and hence the system can be solved for the $g_u=g_u(a)$
and the desired support $U(a)=\{u\in U\,:\,g_u(a)\neq 0\}$ of
$\det(M_1(\gamma)M_2)$ can indeed be computed in polynomial time.
\qed
\end{proof}

Next, we demonstrate that finding an optimal common base for a
nonlinear matroid intersection problem can be reduced to finding an
optimal $W$-image for a small number of subproblems. Consider data
for a nonlinear matroid intersection problem, consisting of two
matroids $M_1$ and $M_2$ on the common ground set
$N=\{1,\ldots,n\}$~, weight matrix $W\in \Z^{d\times n}$~, and
function $f:\R^d\rightarrow\R$~. Each subset $S\subseteq N$ gives a
\emph{subproblem} of nonlinear matroid intersection as follows. The
matroids of the subproblem are the {\em restriction} of $M_1$ and
$M_2$ to $S$~; that is, the matroid $M_i.S$, $i=1,2$, on ground set
$S$ in which a subset $I\subseteq S$ is independent if and only if
it is independent in $M_i$~. Note that the {\em restriction} of the
matrix representing $M_i$ to the columns indexed by $S$ is the
matrix representation of $M_i.S$~. Then, the subproblem is the
nonlinear matroid intersection problem, consisting of two matroids
$M_1.S$ and $M_2.S$ on the common ground set $S$~, the new weight
matrix is the restrictions of the original weight matrix to $S$~,
and the function $f:\R^d\rightarrow\R$ is the same as in the
original problem. We have the following useful statement.

\begin{lemma}\label{Objective} The nonlinear matroid intersection problem of finding
an optimal common base of two matroids $M_1$ and $M_2$ is reducible
in time polynomial in $n$ to finding an optimal $W$-image for at
most $n+1$ subproblems. \end{lemma}
\begin{proof} Denote by $u^*(S)$ the optimal $W$-image for the matroid intersection subproblem
defined by the matroids $M_1.S$ and $M_2.S$ on  the ground set
$S\subseteq N$~. Now compute an optimal common base for the original
problem, by computing an optimal $W$-image for $n+1$ such
subproblems as follows.

\restylealgo{plain}
\begin{algorithm2e}[H] \vspace{0.2cm} Start with
$S:=N$\; Compute an optimal $W$-image $u^*:=u^*(N)$ of the
original problem\; \For{j=1,2,\dots,n}{ Let
$T:=S\setminus\{j\}$\;Compute $rank(M_1.T)$ and $rank(M_2.T)$\;
\If {$rank(M_1.T)=r$ and $rank(M_2.T)=r$ } {Compute an optimal
$W$-image $u^*(T)$\; \lIf {$f(u^*(T))\geq f(u^*)$}{let $S:=T$};
}}{\bf return} $B:=S$\;
\end{algorithm2e}

\noindent It is not hard to verify that the set $B$ obtained is
indeed an optimal common base for the original problem.\qed
\end{proof}

We are now in position to state our main theorem of this section.
By a \emph{randomized algorithm} that solves the nonlinear matroid
intersection problem we mean an algorithm that has access to a
random bit generator and on any input to the problem outputs a
common base that is optimal with probability at least a half. The
running time of the algorithm includes a count of the number of
random bits used. Note that by repeatedly applying such an algorithm
many times and picking the best common base, the probability of
failure can be decreased arbitrarily; in particular, repeating it $n$
times decreases the failure probability to as negligible a fraction
as $1/2^n$ while increasing the running time by a linear factor
only.

\begin{theorem} \label{Random}
For every fixed $d$~, there is a randomized algorithm that, given an
integer weight matrix $W\in\Z^{d\times n}$~, and a function
$f:\R^d\rightarrow\R$ presented by a comparison oracle, solves the
nonlinear matroid intersection problem in time polynomial in $n$ and
$\max W_{i,j}$~.
\end{theorem}

\begin{proof}
As explained at the beginning of this section, we may and will
assume that $W$ is nonnegative. We claim that the following
adaptation of the algorithm of Lemma \ref{Objective} provides a
common base which is optimal with probability at least $1/2$~.
Apply the algorithm of Lemma \ref{Objective}, but at each
iteration apply the algorithm of Lemma \ref{Interpolation} to
determine a random optimal $W$-image $\hat{u}^*(T)$ which we claim
to be the optimal $W$-image $u^*(T)$ with probability at least
$1-1/\left(2(n+1)\right)$~.  At each iteration, the algorithm
either returns: a correct optimal $W$-image with probability at
least $1-1/\left(2(n+1)\right)$~; or a wrong one, or the computed $\hat{U}$
is empty with probability at most
$1/\left(2(n+1)\right)$~. In case that the $\hat{U}$ computed by the
randomized algorithm is empty, we set $\hat{u}^*(T)$ to be a
virtual point with objective value $f(\hat{u}^*(T))=-\infty$~.

Now we show that the probability of success in each iteration,
that is, of computing an optimal $W$-image $u^*(T)$ for any
$T\subseteq N$ in which $rank(M_1.T)=rank(M_2.T)=r$~, is at least
$1-1/\left(2(n+1)\right)$~, as claimed before. Indeed, using
polynomially many random bits, draw independent and uniformly-distributed
integers from $\{1,2,\dots,2r(n+1)\}$ and substitute
them for the $a_{j}$~. Next compute a point $\hat{u}^*(T)\in
\widehat U$ attaining $\max\{f(u):u\in \widehat U\}$ using the
algorithm underlying Lemma \ref{Interpolation} and a comparison
oracle. By Proposition \ref{Probability}, with probability at
least $1-1/\left(2(n+1)\right)$~, we have $u^*(T)\in\widehat U$~, in
which event $\max\{f(u):u\in \widehat U\}=\max\{f(u):u\in U\}$ and
then $u^*(T)$ is indeed an optimal $W$-image.

Therefore, the probability that the algorithm obtains a correct
optimal $W$-image at all iterations, and consequently outputs an
optimal common base $B$ of $M_1$ and $M_2$~, is at least the
product over all iterations of the probability of success in each
iteration, that is $(1-1/\left(2(n+1)\right))^{n+1}$~.  Hence at
least $1/2$ as desired. This completes the proof.

\qed
\end{proof}

\begin{acknowledgements}
Yael Berstein was supported by the Israel Ministry of Science
Scholarship for Women and by a scholarship from the Graduate
School of the Technion. The research of Jon Lee, Shmuel Onn and
Robert Weismantel was partially supported by the Mathematisches
Forschungsinstitut Oberwolfach during a stay within the Research
in Pairs Programme. Shmuel Onn was also supported by the ISF -
Israel Science Foundation. Robert Weismantel was also supported by
the European TMR Network ADONET 504438.
\end{acknowledgements}

\bibliographystyle{spmpsci}

\bibliography{MatroidInt_report}

\end{document}